\documentclass[11pt]{article}
\usepackage{amsmath,amssymb}
\usepackage{color}
\def\comm#1{\textcolor{red}{\marginpar{ZZ}[#1]}}
\RequirePackage{graphics,epsfig,psfrag}
\def\abs#1{\left \vert #1 \right \vert}

\def\RR{{\bf R}} 
\def\ZZ{{\bf Z}} 
\def\cS{{\cal S}} 
\def\Mod#1{\,(\hbox{\rm mod}\,#1)}

\def\TT{{\rm T}}

\def\tr{\hbox{\rm tr}\,}

\def\phi{\varphi}
\def\eps{\varepsilon}

\def\cC{{\cal C}}
\def\card{\hbox{\rm card}\,}

\def\undemi{{1\mskip-3mu /2}}
\def\pn{\medskip\par\noindent}
\def\bi{\vspace{-2pt}\begin{itemize}\itemsep -2pt plus 1pt minus 1pt}
\def\ei{\end{itemize}\vspace{-4pt}}
\def\bn{\vspace{-2pt}\begin{enumerate}\itemsep -2pt plus 1pt minus 1pt}
\def\en{\end{enumerate}\vspace{-4pt}}
\newcommand{\Pf}{{\rm Proof}. }

\newcommand{\EPf}{\hbox{}\hfill$\Box$\vspace{.5cm}}
\def\[#1\]{\begin{eqnarray}#1\end{eqnarray}}
\def\$#1\${\begin{eqnarray*}#1\end{eqnarray*}}

\def\pent#1#2{\lfloor\frac{#1}{#2}\rfloor}

\def\comb#1#2{{}{\textstyle{{#1} \choose {#2}}}}

\def\Sum{\mathop{\sum}\limits}

\def\abs#1{\left \vert #1 \right \vert}

\def\frac#1#2{{\textstyle{{#1} \overwithdelims.. {#2}}}}
\def\Frac#1#2{{\displaystyle{{#1} \overwithdelims.. {#2}}}}

\catcode`\@=11
\def\system#1{\left\{\null\,\vcenter{\openup\jot\m@th
\ialign{
\strut\hfil$\displaystyle{##}$&
$\,\displaystyle{{}##}\,$\hfil&&
\strut\hfil$\,\displaystyle{##}$&
$\,\displaystyle{{}##}\,$\hfill
\hfil\crcr#1\crcr}}\right.}

\def\cmatrix#1{\left [
\null\,\vcenter{
\ialign{
\hfil${##}\ $\hfil &
\hfil$\ {##}\ $\hfil&&
\hfil$\ {##}\ $\hfil&
\hfil$\ {##}$\hfil
\crcr#1\crcr}}\right ]}

\def\@begintheorem#1#2#3{\par\addvspace{8pt plus3pt minus2pt}%
              \noindent{\csname#1headfont\endcsname#1\ \ignorespaces#3 #2.}%
              \csname#1font\endcsname\hskip6pt\ignorespaces}
\def\@endtheorem{\par\addvspace{8pt plus3pt minus2pt}\@endparenv}

\def\blfootnote{\xdef\@thefnmark{}\@footnotetext} 
\catcode`\@=12

\newtheorem{theorem}{Theorem}[section]
\newtheorem{thm*}{Theorem}
\newtheorem{thm}[theorem]{Theorem}
\newtheorem{corollary}[theorem]{Corollary}
\newtheorem{lemma}[theorem]{Lemma}
\newtheorem{proposition}[theorem]{Proposition}

\newtheorem{remark}[theorem]{Remark}
\newtheorem{remarks}[theorem]{Remarks}
\newtheorem{example}[theorem]{Example}

\newtheorem{algorithm}[theorem]{Algorithm}
\newtheorem{conjecture}[theorem]{Conjecture}
\def\[#1\]{\begin{align}#1\end{align}}

\date{\today}
\title{Conway polynomials of two-bridge links}
\author{Pierre-Vincent Koseleff, Daniel Pecker}
\begin{document}
\maketitle
\begin{abstract}
We give necessary conditions for a polynomial to be the Conway
polynomial of a two-bridge link. As a consequence, we obtain simple
proofs of the classical theorems of Murasugi and Hartley. We give a
modulo 2 congruence for links, which implies the classical modulo 2
Murasugi congruence for knots. We also give sharp bounds for the
coefficients of the Conway and Alexander polynomials of a two-bridge
link. These bounds improve and generalize those of Nakanishi and Suketa.
\end{abstract}
\begin{center}
\parbox{12cm}{\small
\tableofcontents }
\end{center}
\section{Introduction}
In this paper, we study the problem of determining whether a given polynomial is the Conway
polynomial of a two-bridge link or knot.
For small degrees, this problem can be solved by an exhaustive search
of possible two-bridge links.
Here, however, we give necessary conditions
on the coefficients of the polynomial, which can be tested for high degree polynomials.
\pn
In section \ref{sec:conway} we present Siebenmann's description of the Conway polynomial of a two-bridge link.
We obtain a characterization of modulo 2 two-bridged Conway polynomials
with the help  of the Fibonacci polynomials $f_k$ defined by:
\[ f_0=0, f_1=1, f_{n+2} (z)= z f_{n+1}(z) + f_n(z), \, n \in \ZZ.\label{fibo}\]
\pn
{\bf Theorem \ref{fib2}.}
{\em
Let $\nabla (z) \in \ZZ[z]$ be the Conway polynomial of a rational  link (or knot).
There exists a Fibonacci polynomial $f_D(z)$ such that $ \nabla(z) \equiv f_D(z) \Mod 2$.}
\pn
We give a simple method (Algorithm \ref{degree}) that determines
this Fibonacci polynomial.
\pn
In section \ref{sec:inequalities}, we obtain inequalities for the coefficients of
the Conway polynomials of links (or knots)  denoted by
$$ \nabla_m (z) = \sum_{k=0}^{\pent m2 } c_{m-2k} z^{m-2k}.$$

\pn
{\bf Theorem \ref{nak1}.}
{\em
For $ k\geq 0$,
$$\abs {c_{m-2k}} \le   \comb{ m-k}{k }   \abs {c_m}.$$
If equality holds for some  positive integer $k< \pent m2$,
then it holds for all integers. In this case,  the link is  isotopic to a link
 of Conway form $ C ( 2, -2, 2,  \ldots, (-1)^{m+1} 2 )$
 or $ \allowbreak C(2,2, \ldots, 2) $, up to mirror symmetry.}
\pn
When $ |c_m| \ne 1,$ we have the following sharper bounds:
\pn
{\bf Theorem \ref{nak3}.}
{\em Let $ g \ge 1$ be the greatest prime divisor of $c_m, $ and let $k \ne 0$. Then
$$
\abs { c_{m-2k}} \le \Bigl( \comb{ m-k-1}{k}  + \Frac 1g \bigl( \comb{ m-k-1}{k-1} -1 \bigr)
\Bigr) \abs {c_m} +1.
$$
Equality holds for links of Conway forms
$C (2g, 2,2, \ldots, 2 )$  and $ C( 2g, -2,2, \ldots, (-1)^{m +1}\, 2)$.}
\pn
We also obtain the following trapezoidal property:
\pn
{\bf Theorem \ref{alpha}.}
{\em Let $K $
be a two-bridge  link (or knot).
Let
$$
\nabla_K= c_m
 \Bigl(\sum_{i=0}^{\pent m2} (-1)^i \alpha_{i} f_{m-2i+1}\Bigr),  \ {}   \
\alpha_0=1
 $$
be its Conway polynomial expressed in the Fibonacci basis. Then we have
\bn
\item
$\alpha_{j}\geq 0, \ j = 0,\ldots,\pent m2$.
\item If $\alpha_{i}=0$ for some $i>0$ then
$\alpha_{j}=0$ for $j\geq i$.
\en
}
\pn
In section \ref{sec:applications}, we apply our results to the Alexander polynomials.
Theorem \ref{fib2} provides  an easy proof of a congruence of Murasugi
\cite{Mu2} for two-bridge knots.
Moreover, we also obtain  a  congruence for the Hosokawa polynomials
of two-bridge links.
\pn
Then, as a consequence of Theorem \ref{alpha}, we obtain a simple proof
of both the Murasugi alternating theorem (\cite{Mu,Mu1}), and the
Hartley trapezoidal theorem (\cite{Ha}, see also \cite{Kan1}).
\pn
We conclude this section by giving bounds for the coefficients of
the Alexander coefficients.
These  bounds improve those  of Nakanishi and Suketa  for  Alexander polynomials of
two-bridge knots (see \cite[Theorems~2 and 3]{NS}).
Moreover, they are sharp and    hold for any $k$.
\pn
We prove that the conditions on the
Conway coefficients are better than the conditions
on the Alexander coefficients deduced from them.
\pn
In section \ref{sec:conjecture}, we conclude our paper
with the following convexity conjecture:
\pn
{\bf Conjecture \ref{con-con}.}
{\em Let $P(t) = a_0 -a_1 ( t+ t^{-1}) + a_2 ( t^2 + t^{-2} ) - \cdots +(-1)^n a_n ( t^n + t^{-n} )$ be
the Alexander polynomial of a two-bridge knot. Then there exists an integer $k \le n$ such that
$ (a_0, \ldots, a_k) $ is convex and $(a_{k}, \ldots, a_n)$ is concave.}
\pn
We have tested this conjecture for all two-bridge knots with 20 crossings or fewer.

\section{Conway  polynomial}\label{sec:conway}
\begin{figure}[th]
\begin{center}
{\scalebox{.8}{\includegraphics{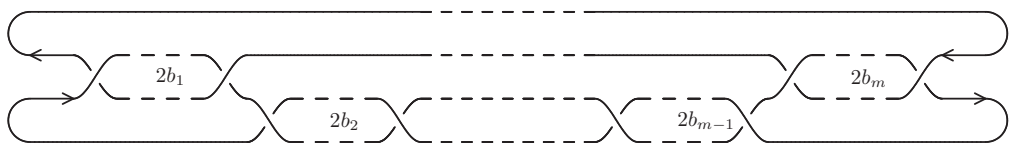}}}
\end{center}
\caption{Oriented two-bridge links ($m$ odd)}
\label{conways3}
\end{figure}

Any oriented two-bridge link can be put in the form shown in Figure \ref{conways3}.
It will be denoted by $C(2b_1, 2b_2, \ldots, 2b_m)$
with $b_i \ne 0$ for all $i$,
including the indicated orientation (see \cite[p.~26]{Kaw},
\cite{Ko,KM}). This is a two-component link if and only if $m$ is odd.

Its Conway polynomial $\nabla _m$ is then given by the Siebenmann
method (see \cite{Si,Cr}).
\psfrag{a}{\small $2b_1$}\psfrag{b}{\small $2b_2$}%
\psfrag{c}{\small $2b_{m-1}$}\psfrag{d}{\small $2b_{m}$}%
\begin{thm}[Siebenmann, \cite{Si,Cr}]\label{th:siebenmann}
Let $ \nabla_m= \nabla_m(z)$ be the Conway polynomial of the oriented two-bridge link
(or knot) of Conway form $ C ( 2b_1, -2b_2, \ldots, (-1)^{m+1} 2 b_m)$.

Let $ \nabla _ {-1}=0, \   \nabla_0 =1$.Then
$
\nabla_ m= b_mz \nabla_{m-1} + \nabla_{m-2}
$
for $m \ge 1$.
\end{thm}
When $z=1,$ this is the  classical Euler continuant polynomial (see \cite{Kn}).
When all the $b_i$ are equal to 1, we obtain the Fibonacci polynomials.

\begin{example}  {\bf The torus links $\TT(2,m)$}.
The Conway polynomial of the torus link $\TT(2,m)= C (2, -2, \ldots , (-1)^{m+1} 2) $ is the Fibonacci polynomial
$f_{m}(z)$ (see  \cite{Ka,KP4}).
\end{example}
Consequently, the following result gives in fact a characterization of
modulo 2 Conway polynomials of two-bridge links.

\begin{thm}{\label{fib2}}
Let $\nabla_m$ be the Conway polynomial of a two-bridge link. 
Then there exists a Fibonacci polynomial $f_D$ such that $\nabla_m \equiv f_D \Mod 2$.
\end{thm}
\Pf
Let us  write
$ (a,b) \equiv (c,d) \Mod 2 $ when $ a \equiv c \Mod 2$ and $b \equiv d \Mod 2$.
We will show by induction on $m$ that there exist integers $D$ and $e = \pm 1$ such that
$(\nabla_{m-1},\nabla_{m}) \equiv (f_{D-e},f_{D}) \Mod 2$.
\pn
The result is true for $m=0$ as $(\nabla_{-1},\nabla_0) = (0,1) = (f_0,f_1)$, that is $D=e=1$.
\pn
Suppose that  $(\nabla_{m-1},\nabla_{m}) \equiv (f_{D-e},f_{D}) \Mod 2$, with $e=\pm 1$
for some $m\geq 0$.
Then we have $\nabla_{m+1} = b_{m+1} z \nabla_{m} + \nabla_{m-1}$.
\par
If $b_{m+1}\equiv 0 \Mod 2$ then
$
\nabla_{m+1} \equiv  \nabla_{m-1} \equiv  f_{D-e} \Mod 2
$
and $(\nabla_{m},\nabla_{m+1}) \equiv (f_{D},f_{D-e})$.
\par
If $b_{m+1}\equiv 1 \Mod 2$ then
$
\nabla_{m+1} \equiv  z f_{D} + f_{D-e} \equiv f_{D+e} \Mod 2,$
and  consequently $(\nabla_{m},\nabla_{m+1}) \equiv (f_{D},f_{D+e})$.
\EPf
\begin{example}
The Pretzel knot $8_5$ has Conway polynomial $1-z²-3z^4 - z^6 \equiv
f_1 + f_3+f_7 \Mod 2$.  By theorem \ref{fib2}  it is not a two-bridge knot.
\end{example}
From the proof of Theorem \ref{fib2}, we deduce a fast algorithm
for the determination
of the integer $D$ such that $\nabla_K \equiv f_D \Mod 2$, see also \cite{Bu}.
\begin{algorithm}\label{degree}
Let $K$ be a two-bridge link (or knot) of Conway form $C( 2b_1, 2b_2, \ldots, 2b_m)$.
Let us define the sequences of integers $e_i$ and $D_i$, $ i=0,\ldots, m,$ by
$$e_0=1, \  D_0 = 1, \ e_{i+1} = -(-1)^{b_{i+1}}e_i, \ D_{i+1} = D_i + e_{i+1}.$$
Then we have  $ \nabla (z) \equiv   f_D (z) \Mod 2$ where $D= \abs {D_m}$.
\end{algorithm}
This algorithm may be useful for the study of Lissajous knots.  Jones, Przytycki, and Lamm proved
that the Conway polynomial of a  two-bridge Lissajous knot satisfies
the congruence
$ \nabla (z) \equiv 1 \Mod2,$ that is $D=0$ (see \cite{BDHZ,JP,La}).

\section{Inequalities for the coefficients of the Conway
  polynomial}\label{sec:inequalities}
We shall need the following  explicit notation for Conway polynomials:
$$
\nabla_m (z)= \sum _{k=0}^{ \pent m2} c_{m-2k} (b_1, \ldots, b_m) z^{m-2k}.
$$
Thus, the Siebenmann formula (Theorem \ref{th:siebenmann})
means that
\[
c_{m-2k} (b_1, \ldots , b_m)=
b_m \cdot c_{m-1-2k}(b_1, \ldots , b_{m-1} )
+ c_{m-2k} ( b_1, \ldots , b_{m-2}). \label{continu}
\]

\begin{remark}\label{fibex}
For the torus link $ T(2,m)= C (2,-2, \ldots, (-1)^{m+1} 2),$ all the $b_i$ are equal to $1$,
 and an easy induction shows that $ c_{m-2k}(1, \ldots , 1 )= \comb {m-k}k.$
consequently,
we obtain the following expression for the Fibonacci polynomials:
$$
f_{m+1}(z) = \sum_{k=0}^{ \pent m2}
\comb{ m-k}{k} z ^{m-2k} \quad {\rm for} \quad m \ge 0  .
$$
This means that the Fibonacci polynomials can be read on the
diagonals of  the \comm{the}
Pascal's triangle.
When $z=1$, we recover the classical Lucas identity
$$
F_m = \sum_{k=0}^{ \pent m2} \comb{ m-k}{k},
$$
where $F_m$ are  the Fibonacci numbers ($F_0=0, \  F_1=1,  F_{n+1}=F_n + F_{n-1}$).
\end{remark}
In the next result, we deduce some properties of the coefficient
$ c_{m-2k} (b_1, \ldots, b_m),$
considered as a polynomial in the $m$ variables $b_1, \ldots,  b_m.$
\begin{proposition}\label{euler}
Let $\cC(m,k),$  $ \  m \ge 2k,$ be the set of all monomials
$\Frac {b_1 \cdots b_m}{ b_{i_1}b_{i_1+1} \cdots b_{i_k} b_{i_k +1}},$
where
$i_h +1 < i_{h+1}$. Let $\cC_j(m,k)$ be the subset of all monomials
of $\cC(m,k)$ that are relatively prime to $b_j.$
Then we have
\bn
\item  The set $\cC (m,k)$ has
 $\comb{ m-k}{k}$ elements.
\item  The polynomial
$c_{m-2k} (b_1, \ldots, b_m) $ is the sum of all monomials of $\cC (m,k)$.
\item If $k \ne 0$, then the monomials of $\cC (m,k)$ do not have a common divisor except 1.
\item The number of elements of $\cC_j(m,k)$ is at least
 $ \comb{ m-1-k}{k-1}$.
\item If $k \ge 2$, then the monomials of $\cC_j (m,k)$
do not have a common divisor except 1.
\en
\end{proposition}
\Pf
\bn
\item
By induction on $m.$

We have
$\cC(m,0)= \{ \, b_1\cdots b_m  \,  \}$, $\cC(2,1) = \{ 1\}$ and
$\cC(3,1)= \{ b_1, \, b_3  \}$.
Hence the result is true for $k=0,$ and also for $m \le 3.$

Let us suppose the result true for $m-1$ and $m-2.$
We can suppose $k \ne 0.$
If a monomial of $ \cC(m,k)$ is not a multiple of $b_m$, then it is
not a multiple of $b_{m-1}  $ either, and consequently it is an
element of $ \cC (m-2,k-1).$ 
Therefore, we have the following partition of $ \cC(m,k)$ for $k\ne 0$:
$$
\cC(m,k) = b_m \cdot \cC(m-1,k) \ \bigsqcup \ \cC(m-2,k-1),
$$
and then
$$
\card \cC(m,k)  =  \card \cC(m-1,k)
+ \card \cC(m-2,k-1) =
\comb {m-1-k}k + \comb {m-1-k} {k-1} = \comb {m-k}k.
$$
\item
By induction on $m.$ Using our partition of $\cC(m,k),$ we see that the sum of
the monomials of $\cC(m,k)$ is
$b_m \, \cdot \, c_{m-1-2k} ( b_1, \ldots , b_{m-1})  +
 c_{m-2k} (b_1, \ldots , b_{m-2}).$

By Siebenmann's formula, this polynomial is equal to $ c_{m-2k} (b_1, \ldots , b_m).$
\item If $k \ne 0,$ then for every integer $i \le m,$ there is an element of $\cC(m,k)$
which is not divisible by  $b_i.$
Hence the GCD of the elements of $\cC(m,k)$ is 1.
\item
Let
$\mathbf{b} = (1,\ldots,1,0,1,\ldots,1) \in \RR^m$ where $b_j=0$, and
$b_k=1$ for $k\ne j$.

Let us define the polynomials $g_n, $ for $n \le m$ by
$g_n(z)= \nabla_n (\mathbf{b}) (z)$.
The number of  elements of $\cC_j(m,k)$  is the coefficient
$ c_{m-2k}(\mathbf{b})$ of $g_m(z)$.
\pn
If $j=1$, then we have $g_1=0$, $g_2=1,$ and $ g_n=z g_{n-1}+ g_{n-2}$ for $ n \ge 2.$
Then, an easy induction shows that
 $g_n=f_{n-1}. $
\pn
If $j>1$, then we have
$g_1 = f_2, \  \ldots, \  g_{j-1}= f_{j}, \   g_j = f_{j-1}$,
and
$g_{n+1}= z g_n + g_{n-1}$ for $ \ n \ge j$.

Let us  write $ p(z)\succeq q(z)$ when each coefficient of $p$ is
greater than or equal to the corresponding coefficient of $q$.
We have $f_{k+2} \succeq f_k$, and therefore
$ g_{j+1}= z f_{j-1} + f_j \succeq z f_{j-1} + f_{j-2}= f_j .$
Then a simple induction shows that $g_m \succeq f_{m-1}$, and consequently
$ c_{m-2k}( {\bf b}) \ge \comb {m-1-k}{k-1}.$
\item
Since $ k  \ge 2,$  then for every  $i \ne j$ there is
a monomial of $\cC_j (m,k)$ which is not divisible by  $b_i$.
Consequently, the GCD of the elements of
 $\cC_j (m,k)$  is 1.
\EPf
\en

\begin{thm}\label{nak1}
For $ k\geq 0$,
$$\abs {c_{m-2k}} \le   \comb{ m-k}{k }   \abs {c_m}.$$
If equality holds for some
integer $k< \pent m2$,
then it holds for all integers. In this case,  the link is  isotopic
to the torus link $T(2,m) $
or to the link $ C(2,2, \ldots, 2) $, up to mirror symmetry.
\end{thm}
\Pf
By Proposition \ref{euler},  the number of monomials of
$c_{m-2k} (b_1, \ldots, b_m)$  is
$ \comb{ m-k}{k}$. The result follows since no monomial
 is greater than $\abs{c_m}= \abs{ b_1\cdots b_m}$.
\pn
If equality holds for some positive integer $k< \pent m2$,  then for
all $i, \, j, \  $ $ b_i b_{i+1}= b_j b_{j+1} = \pm 1,$ which implies
the result.
\EPf
\begin{example}
The knot $10_{145}$ has Conway polynomial $P=1+5z^2+z^4$. We have $P
\equiv f_5 \Mod 2,$ but $P$ does not satisfy the condition $\abs{c_2}
\leq 3$, and then $10_{145}$ is not a two-bridge knot.

The knot $11n109$ has Conway polynomial $1+6z^2+z^4 - z^6 .$
It satisfies the bounds of Theorem \ref{nak1}: $\abs{c_2} \leq 6, \
\abs{c_4} \leq 5$, but not the equality condition:  $c_2=6$ whereas $c_4 \ne 5$.
Consequently, $11n109$ is not a two-bridge knot.
\end{example}
To prove the refined inequalities of Theorem \ref{nak3}, we shall use
the following lemma, which generalizes the inequality
$a+b \le ab+1$,  valid for positive integers (see also \cite{NS}).
\begin{lemma}\label{fact1}
Let $p_i({\bf  x}), i \in \cS $  be relatively prime  divisors of
$p({\bf x}) = x_1 x_2 \cdots x_m$.

Let  $ {\bf b}= (b_1, \ldots, b_m) $ be a m-tuple of positive integers.
Then
\[
\sum _{i \in \cS} p_i ( {\bf b})  \le
\Bigl( \card (\cS)-1 \Bigr) {p( {\bf b})} +1.\label{fact2}
\]
\end{lemma}
\Pf
We do not suppose that the $p_i$ are distinct integers.
Let us prove the result by induction on $k={\rm card}(\cS)$.
If $k=1,$ then we have ${p_1}=\pm 1,$ and
the result is true.
When  all the $p_i$ are equal to 1,  the result is true.
Otherwise, let $x_h$ be a divisor of some $p_i$.
\pn
Let $\cS_1= \{ i\in \cS : \   x_h \mid p_i \} $ and $\cS_2= \cS-\cS_1$.
We have $k=k_1+k_2$,
where $k_j= {\rm card} (\cS_j)$.
Let $q_j= {\rm GCD}\{ p_i, i \in \cS_j \}$, then $q_1$ and $q_2$ are coprime,
and $q_1q_2$ is a divisor of $p$.
\pn
By induction we obtain for $j= 1, 2$:
$$ \sum _{i \in \cS_j} p_i( {\bf b}) \le {q_j( {\bf b})}
\Bigl( ( k_j-1) \Frac { {p( {\bf b} )}} {{q_j ( {\bf b})} } +1 \Bigr)
= (k_j-1)p ({\bf b}) + q_j( {\bf b }) .$$
Adding these two inequalities we get
\begin{align*}
\sum _{ i \in \cS} p_i ( {\bf b}) &\le
( k_1+k_2-1) {p( {\bf b})} + {q_1 ( {\bf b})} + {q_2 ( {\bf b})} - {p({\bf b})}\\
&\le
(k-1) {p({\bf b})} + {q_1({\bf b})} {q_2( {\bf b})}  - {p({\bf b})}+1,
\end{align*}
which proves the result, since ${q_1({\bf b})}{ q_2( {\bf b})} \le {p({\bf b})}$.
\EPf
\pn
With this lemma we can prove:
\begin{thm}\label{nak3}
Let $ g \ge 1$ be the greatest prime divisor of $c_m$, and let $k \ne 0.$ Then
$$
\abs { c_{m-2k}} \le \Bigl( \comb{ m-k-1}{k}
 + \frac 1g \bigl( \comb{ m-k-1}{k-1} -1 \bigr)
\Bigr) \abs {c_m} +1.
$$
Equality holds for links of Conway form
$C(2g,-2, \ldots, (-1)^{m+1} 2)$ and $C(2g,2, \ldots,  2)$.
\end{thm}
\Pf
If $k=1$, then by Proposition \ref{euler}
 the polynomial
$c_{m-2}(b_1,\ldots,b_m)$ is the sum of $m-1$ coprime monomials.
Then, using Lemma \ref{fact1} and the notation
$ \abs {\bf b} =(\abs {b_1}, \ldots, \abs {b_m} )$, we get
$$
\abs{c_{m-2}} = \abs {c_{m-2}({\bf b})} \le c_{m-2}( \abs {\bf b} )
 \leq (m-2) c_{m}(\abs{\bf b}) + 1 = (m-2) \abs{c_m} + 1.
$$
Now, suppose $ k \ge 2$. Let $g$ be the greatest prime divisor of
the integer $c_m= b_1 \cdots b_m,$ and suppose that $g \mid b_j$.
Let $ {\cal M} $ be the set of monomials of
 $c_{m-2k}( b_1, \ldots, b_m), $ and let
 $ {\cal M}_j$ be the subset of monomials of ${\cal M}$
that are prime to  $b_j$.

By Proposition \ref{euler},
the monomials of ${\cal M}_j$ are relatively prime, and their number $N$ verifies
$N \ge \comb{m-1-k}{k-1}$. Using Lemma \ref{fact1} we obtain:
$ \sum_{ p_i \in
{\cal M}_j} p_i ({\bf b}) \le (N-1) \Frac{\abs{ c_m} }{\abs{b_j}} +1$
and then
\begin{align*}
\abs {c_{m-2k}} = \abs {\textstyle{\sum _{ p_i \in {\cal M}}} p_i
({\bf b})}
&\le \Bigl( \Frac{N-1}{g}+ (\comb{m-k}{k}-N ) \Bigr)\abs {c_m}+1\\
&=
\Bigl( \comb{m-k}{k} -N (1 - \frac 1g) - \frac 1g\Bigr)\abs {c_m}+1\\
&\le
\Bigl( \comb{ m-k}{k} - \comb{ m-1-k}{k-1} (1-\frac 1g) - \frac 1g  \Bigr)\abs {c_m}+1\\
&=\Bigl( \comb{ m-1-k}{k} + \frac 1g  (\comb{ m-1-k}{k-1}-1 ) \Bigr)
\abs {c_m} +1.
\end{align*}
For links of Conway form $ C(2g, 2, \ldots, 2 )$
or $C(2g,-2, \ldots, (-1)^{m+1} 2)$, we have
$ \abs {\mathbf{b} } = (g,1,\ldots,1)$, $N=\comb{ m-1-k}{k-1}$, $\abs {c_m}=g,$
$\abs {c_{m-2k} }= g \comb {m-1-k}{k} + \comb {m-1-k}{k-1},$ and equality holds.
\EPf
\pn
We will now express the Conway polynomials of two-bridge links in
terms of Fibonacci polynomials, and show that their coefficients
are alternating.
\begin{thm}\label{alpha}%
Let $K$ be a two-bridge  link (or knot).
Let
$$
\nabla_K=
c_m
 \Bigl(\sum_{i=0}^{\pent m2} (-1)^i \alpha_{i} f_{m-2i+1}\Bigr),  \ {}   \
\alpha_0=1
 $$
be its Conway polynomial written in the Fibonacci basis. Then we have
\bn
\item
$\alpha_{j}\geq 0, \ j = 0,\ldots,\pent m2$.
\item If $\alpha_{i}=0$ for some $i>0$ then
$\alpha_{j}=0$ for $j\geq i$.
\en
\end{thm}
\Pf
Let $K= C(2b_1, -2b_2, \ldots, (-1)^{m+1} \, 2b_m), $ with $b_i\ne 0$ for all $i,$
and let $\nabla_n$ be the polynomials obtained in the Siebenmann method.

We have $\nabla_0 = f_1$, $\nabla_1 = b_1 f_2$, $\nabla_2 = b_1b_2
\Bigl(f_3 - (1-\frac{1}{b_1b_2}) f_1\Bigr)$.

Let us show by induction that if
$$
\nabla_m = b_1 \cdots b_m \Bigl(\Sum_{i=0}^{\pent m2} (-1)^i \alpha_{i} f_{m+1-2i}\Bigr), \,
\nabla_{m-1} = b_1 \cdots b_{m-1} \Bigl(\Sum_{i=0}^{\pent{m-1}2}
(-1)^i \beta_{i} f_{m-2i}\Bigr)$$
then $\alpha_{j} \geq  \beta_{j}\geq 0$, and if $\alpha_i=0$ for some $i,$ then
$\alpha_j=0$ for $j \ge i$.
\pn
The result is true for $m=1$ and for $m=2.$
Using $zf_{m+1-2i} =  f_{m+2-2i}-f_{m-2i}$ and $\nabla_{m+1} =
b_{m+1}z \nabla_m + \nabla_{m-1}$,
we deduce that
$$\nabla_{m+1} = b_1 \cdots b_{m+1} \Bigl(\Sum_{i=0}^{\pent {m+1}2}
(-1)^i \gamma_{i} f_{m+2-2i}\Bigr),$$ where
$\gamma_{0} = 1$ and
\[
\gamma_{i}
&= \alpha_{i} + (\alpha_{i-1} - \beta_{i-1}) + (1- \frac{1}{b_m b_{m+1}}) \beta_{i-1},
\, i = 1,\ldots, \pent{m+1}2. \label{abc}
\]
As $\abs{b_m b_{m+1}}\geq 1$, we deduce by induction that $\gamma_{i}
\geq \alpha_{i} \geq 0$.

Furthermore, if $\gamma_i=0,$ then by Formula (\ref{abc}) $
\alpha_i=0,$ and then, by induction, $\alpha_j = \beta_j=0$ for $j\ge
i$.   Finally, by Formula (\ref{abc}), we get $\gamma_j=0$ for $j \ge i$.
\EPf
\section{Applications to the Alexander polynomial}\label{sec:applications}
In this section, we will see that our necessary conditions on the
Conway coefficients imply similar necessary conditions on the
Alexander coefficients of two-bridge knots and links. These
conditions are improvements of the classical results.
\pn
The Conway and the Alexander polynomials of a knot $K$ will be denoted by
$$\nabla_K (z) = 1 +{\tilde c}_1 z^2 + \cdots + {\tilde c}_n z^{2n}$$
and
$$ \Delta_K(t) = a_0 - a_1 (t+t^{-1}) + \cdots + (-1)^n a_n ( t^n+ t^{-n} ).$$
The Alexander polynomial $ \Delta_K (t)$ is deduced from the Conway
polynomial by:
$$ \Delta_K (t) =  \pm \nabla_K \Bigl( t^{\undemi} - t^{ - \undemi} \Bigr). $$
It is often normalized so that  $a_n$ is positive.
Thanks to this formula, it is not difficult to deduce
the Alexander polynomial from the Conway polynomial.
If we use the Fibonacci basis, it is even easier to deduce the
Conway polynomial of a knot from its Alexander polynomial.
\begin{lemma}\label{undemi}
If $ \ z =t^{\undemi} - t^{ -\undemi}, \ $  and $ \  n \in \mathbf{Z} \ $,
then we have the identity
$$  f_{n+1} (z) + f_{n-1} (z) = (t^{\undemi})^n + ( -t^{ - \undemi} )^n,$$
where the $f_k (z)$ are the  Fibonacci polynomials.
\end{lemma}
\Pf
Let
$ A= \cmatrix { z & 1 \cr 1 & 0 } $ be the  (polynomial) Fibonacci matrix.
If   $z =t^{\undemi} - t^{-\undemi}$, then the eigenvalues of $A$
are  $ \  t^{\undemi}$ and $ \  -t^{ -\undemi}$, and consequently
$ \tr A^n =(t^{\undemi})^n + ( -t^{ -\undemi} )^n$.
On the other hand,
we have $ A^n= \cmatrix { f_{n+1}(z) & f_n(z) \cr f_n(z) & f_{n-1}(z) }$, and then
$ \tr A^n = f_{n+1} (z) + f_{n-1} (z)$.
\EPf
\pn
From Lemma \ref{undemi}, we immediately deduce:
\begin{proposition} \label{chvar}
Let the Laurent polynomial $P(t)$ be defined by
$$P(t) = a_0 -a_1 (t+ t^{-1}) + a_2 ( t^2 + t^{-2} ) - \cdots +(-1)^n a_n ( t^n + t^{-n} ).$$
We have
$$P(t)= \sum_{k=0}^n (-1)^k (a_k-a_{k+1}) f_{2k+1}(z),$$
where $z =t^{\undemi} - t^{-\undemi} $ and $a_{n+1}=0$.
\end{proposition}
Using the  substitution $a_0 = \ldots = a_n=1,$
We deduce the following  useful formula.
\[
f_{2n+1} \bigl( t^{\undemi} - t^{- \undemi} \bigr)=
(t^n + t ^{-n}) - (t^{n-1} + t ^{1-n}) + \cdots + (-1)^n.
\label{afib}
\]
Then, we deduce a simple proof of an elegant criterion due to
Murasugi (\cite{Mu2,Bu})
\begin{corollary}[Murasugi (1971)] {\label{al1}}
Let
$\Delta (t) = a_0 -a_1 ( t+ t^{-1}) + a_2 ( t^2 + t^{-2} ) - \cdots
+(-1)^n a_n ( t^n + t^{-n} )$ be the Alexander polynomial of a
two-bridge knot. There exists an integer $k \le n$ such that
$ a_0, a_1, \ldots, a_k$ are odd, and $a_{k+1}, \ldots, a_n $ are even.
\end{corollary}
\Pf
If $K$ is a two-bridge knot, its Conway polynomial is a modulo 2
Fibonacci polynomial $f_{2k+1}$. By Proposition \ref{chvar} we have
$
 f_{2k+1} \bigl( t^{\undemi} - t^{- \undemi} \bigr)=
( t^k + t^{-k} ) - ( t^{k-1} + t^{1-k} ) + \cdots + (-1)^{k},
$
and the result follows.
\EPf
\begin{remark}
This congruence may be used as a simple criterion to prove that some knots
 cannot be two-bridge knots.
There is a more efficient criterion by Kanenobu \cite{Kan2,St00} using
the Jones and Q polynomials.
\end{remark}
\pn
We also deduce an analogous result for two-component links
\begin{corollary}[Modulo 2 Hosokawa polynomials of two-bridge links]
Let $\Delta (t)= \bigl( t^{\undemi} - t^{- \undemi} \bigr) \Bigl(
 a_0 -a_1 ( t+ t^{-1}) + a_2 ( t^2 + t^{-2} ) - \cdots +(-1)^n a_n (
 t^n + t^{-n} ) \Bigr)$
be the Alexander polynomial of a two-component two-bridge link. Then
all the coefficients $a_i$ are even or there exists an integer $k\le
n$ such that $ a_k, a_{k-2}, a_{k-4}, \ldots$ are odd, and the other
coefficients are even.
\end{corollary}
\Pf
If $K$ is a two-component two-bridge link, its Conway polynomial is an odd Fibonacci
polynomial modulo 2, that is of the form $ f_{2h} (z)$.
An easy induction shows that
$$ f_{4k} \bigl( t^{\undemi} - t^{- \undemi} \bigr)=
\bigl( t^{\undemi} - t^{- \undemi} \bigr) \bigl( u_1+u_3+ \cdots + u_{2k-1} \bigr)$$
and
$$ f_{4k+2} \bigl( t^{\undemi} - t^{- \undemi} \bigr)=
\bigl( t^{\undemi} - t^{- \undemi} \bigr) \bigl( 1+ u_2+ \cdots + u_{2k} \bigr),$$
where $u_j= t^j + t^{-j},$ and the result follows.
\EPf
\begin{remark}
This rectifies Satz 4 in \cite[p.~186]{Bu}.
\end{remark}
\pn
Now, we shall show that  Theorem \ref{alpha}  implies both  Murasugi
and  Hartley theorems for two-bridge knots:
\begin{thm}[Murasugi (1958), Hartley (1979)] {\label{al2}}
Let $$P(t) = a_0 -a_1 ( t+ t^{-1}) + a_2 ( t^2 + t^{-2} ) - \cdots +(-1)^n a_n ( t^n + t^{-n} )
, \  a_n>0 $$ be
the Alexander polynomial of a two-bridge knot. There exists an integer $k \le n$ such that
$ a_0=a_1=\ldots =a_k > a_{k+1}> \ldots > a_n$.
\end{thm}
\Pf
Let $K$ be a two-bridge knot and
$\nabla(z)= \alpha_0 f_1 - \alpha _1 f_3 + \cdots + (-1)^n \alpha_{n} f_{2n+1}$
 be its Conway polynomial expressed in the Fibonacci basis.
 By Theorem \ref{alpha}  $\alpha _n \alpha _k \ge 0$ for all $k$, and  if $\alpha _i=0$ for some
 $i$ then $\alpha_j=0$ for $j\le i$.

Let
$\Delta (t) = a_0 -a_1 ( t+ t^{-1}) + a_2 ( t^2 + t^{-2} ) - \cdots +(-1)^n a_n ( t^n + t^{-n} )
, \   a_n>0$
be the Alexander polynomial of $K$.
We have $ \Delta (t)= \eps \nabla ( t^{ \undemi} - t^{- \undemi} ) $,  where $\eps= \pm 1,$
and then, by  Corollary \ref{chvar}, $  \eps \alpha_k = a_k- a_{k+1}$.
We deduce that $ \eps \alpha_n= a_n >0,$ and then
$ a_k- a_{k+1} = \eps \alpha_k \ge 0 $ for all $k$.
Consequently we obtain $ a_0  \ge a_1 \ge \ldots  \ge  a_n >0$.

Furthermore, if $a_k= a_{k-1}$ for some $k,$ then $\alpha_{k-1}=0,$ and consequently
$ \alpha_{j-1}=0$ for all $j\le k$.
This implies that for all $j \le k, $  $a_j=a_{j-1} $, which concludes the proof.
\EPf
\pn
Now, we shall give explicit formulas for  Alexander coefficients in terms of
 Conway coefficients.
\begin{proposition}\label{c2a}
Let $ Q(z)={\tilde c}_{0} + {\tilde c}_{1} z^2 + \cdots + {\tilde c}_{n} z^{2n}$ be a polynomial.
 We have
$$
Q(t^{\undemi} - t^{- \undemi}) =
 a_0 -a_1 ( t+ t^{-1}) + a_2 ( t^2 + t^{-2} ) - \cdots +(-1)^n a_n ( t^n + t^{-n} ),
$$
where
\[
a_{n-j} =  \sum _{k=0}^j (-1)^{n-k}{\tilde c}_{n-k} \comb{2n-2k}{j-k}.\label{c2af}
\]
\end{proposition}
\Pf
It is sufficient to prove Formula (\ref{c2af}) for the monomials
$Q(z)=z^{2m}$.
Let us consider $u_i=t^i + t^{-i}. $ By the binomial formula we have
$$
 \Bigl( t^{\undemi} - t^{- \undemi} \Bigr)^{2m}=
\sum _{k=0} ^{m-1} (-1)^k \comb{2m}{k} u_{m-k} +
(-1)^m \comb{2m}{m}.
$$
and then
$a_{n-j}= (-1)^m \comb {2m}h $ where $ m-h=n-j.$
On the other hand, the proposed formula asserts
$$
a_{n-j} =  \sum _{k=0}^j (-1)^{n-k}{\tilde c}_{n-k} \comb{2n-2k}{j-k} =
(-1)^m \comb {2m}h \quad {\rm where} \quad h=m+j-n,
$$
which is the same result.
%
\EPf
\begin{remark}
Considering the Fibonacci polynomials $f_{2n+1} = \sum_{k=0}^n \comb{2n-k}{k} z^{2n-2k}$,
Formulas (\ref{afib}) and (\ref{c2af}) give
the identity
$$
\sum_{k=0}^j (-1)^{k} \comb{2n-k}{k} \comb{2n-2k}{j-k} \, = \, 1 , \quad  n,j \geq 0  .
$$
\end{remark}
\begin{remark}\label{ra2cf}
Fukuhara  \cite{Fu} gives a converse  formula for the $c_k$ in terms of the $a_k$,
$$
{\tilde c}_{n-j} = \sum _{k=0} ^j (-1)^{n-k} a_{n-k} \frac{2n-2k}{2n-j-k} \comb{2n-j-k}{2n-2j}.
$$
We shall not use this formula.
\end{remark}
From the bounds we obtained for  Conway coefficients we
can deduce a simple proof  of the Nakanishi--Suketa bounds
(\cite[Th. 1, 2]{NS})
for  the Alexander  coefficients.
\comm{}
\begin{corollary}[Nakanishi--Suketa (1993)]\label{naks}
We have the following sharp inequalities
(where all the $a_i$ are positive):
\bn
\item
$a_{n-j} \le a_n \Bigl(  \sum _{k=0} ^j \comb{2n-2k}{j-k}\comb{ 2n-k}{k} \Bigr)$.
\item
$ 2 a_n -1 \le a_{n-1} \le (4n-2) a_n +1$.
\en
\end{corollary}
\Pf
\bn
\item
Using Formula (\ref{c2af}) and Theorem \ref{nak1}, we obtain
\[
\abs{a_{n-j}} \le
\sum_ {k=0}^j \abs{\tilde c_{n-k} } \comb{2n-2k}{j-k} \le
\abs{a_n} \sum_{k=0}^j \comb{2n-k}{k}\comb{2n-2k}{j-k}.\label{ns1}
\]
\item
We have
$ \abs {\tilde c_{n-1}} \le \comb{2n-2}1 \abs { \tilde c_n}  +1$ by Theorem \ref{nak3},
and $a_{n-1} = \tilde c_{n-1} - \comb{2n}1 \tilde c_n$ by Proposition \ref{c2a}. We thus deduce
\[
\abs {a_{n-1}} \le \comb{2n}1 \abs {\tilde c_n} + \comb{2n-2}1 \abs{ \tilde c_n} +1=
(4n-2) \abs{ a_n} +1.\label{ns1b}
\]
We also have
$$
\abs{ a_{n-1}} \ge \comb{2n}1 \abs{\tilde c_{n}} - \abs{ \tilde c_{n-1}} \ge
 \comb{2n}1 \abs{\tilde c_n} - \comb{2n-2}1 \abs{\tilde c_n} -1 = 2 \abs{a_n} -1.
$$

\en
The upper bounds (\ref{ns1}) and (\ref{ns1b}) are attained by the knots $ C(2,2, \ldots, 2)$.
\EPf
\pn
We also have the following sharp bound, which improves
the Nakanishi--Suketa third bound (\cite[Th. 3]{NS})
\begin{thm}\label{nakss}
If $a_n \not = 1,$ then
$a_{n-2} \le ( 8n^2-15n+8) a_n + 2n-1$. This bound is sharp.
\end{thm}
\Pf
From Proposition \ref{c2a} and Theorem \ref{nak3}, we get
\begin{align*}
\abs{a_{n-2}} &\le
\comb {2n}2 \abs{ \tilde c_n} +
\comb {2n-2} 1 \abs{ \tilde c_{n-1}} +
\comb {2n-4}0 \abs {\tilde c_{n-2} }\\
&\le \comb {2n}2 \abs{ \tilde c_n}
+ \comb{2n-2}1 (\comb{2n-2}1 \abs { \tilde c_n}  +1)
+ \Bigl(\comb{2n-3}{2} + \frac 1g ( \comb{2n-3}{1} -1)\Bigr)\abs { \tilde c_n}  +1\\
&=
(8n^2-16n+10 + \frac{2(n-2)}g    ) \abs{a_n} + 2n-1 .
\end{align*}
If $a_n \ne 1$ then $g \ge 2,$ and we obtain
\[
\abs{ a_{n-2}} \le \abs{ a_n} ( 8n^2 -15n +8) + 2n-1. \label{ns2}
\]
This bound  is attained for the knot $C(4,2,2,2, \ldots, 2)$.
\EPf
\pn
The following example shows that the bounds on the Conway coefficients are better than the
bounds on the Alexander coefficients.
\begin{example}
Let us consider the Conway polynomial $ \nabla_ K (z)= 1 +8 z^2 +3 z
^4 -z^6$ of the knot $K=13n1862$ (see \cite{KA}). 
It does not verify the bound of theorem \ref{nak1}, and then it is not a two-bridge knot.
Nevertheless, its Alexander polynomial
$ \Delta _K (t)= 23 - 19 (t+1/t) +9 (t^2+1/t^2) - (t^3+1/t^3) $
 satisfies the bounds of Nakanishi and Suketa,   and also the conditions
of Murasugi and Hartley.
This example shows that the conditions on the Conway coefficients are stronger than
 the conditions on the Alexander coefficient deduced from them.
\end{example}
\begin{remarks}
\bn
\item If $g\geq 3$, we obtained an improvement of the inequality (\ref{ns2}):
$$
a_{n-2} \leq (8n^2-16n+10+ \frac{2(n-2)}g  ) a_n + 2n-1.
$$
\item For $j=3$ we obtain
\begin{align*}
a_{n-3}  &\le
2/3\, \left( 2\,n-3 \right)  \left( 8\,{n}^{2}-24\,n+25 \right) a_n+{
\frac { \left( 3\,n-5 \right)  \left( 2\,n-5 \right)}{g}  a_n}+n \left( 2
\,n-3 \right)\\
&\le
1/6\, \left( 64\,{n}^{3}-270\,{n}^{2}+413\,n-225 \right) a_n+n \left( 2
\,n-3 \right).
\end{align*}
\item Since the inequalities on  Conway coefficients are simpler and stronger,
we shall not give the inequalities on Alexander coefficients for $j \ge 4$.
\en
\end{remarks}

\section{A conjecture}\label{sec:conjecture}
We have computed the Conway polynomials of the $131\,839$ two-bridge
links and knots with 20 or fewer crossings, using Siebenmann's
method. We observed the following property:
\begin{conjecture}\label{trapezoidal}
Let
$\nabla_m = c_m \Bigl(\sum_{i=0}^{\pent m2} (-1)^i \alpha_{i}
f_{m+1-2i}\Bigr), \ \alpha_0=1,$ be the Conway polynomial of a
two-bridge link (or knot) written in the Fibonacci basis. 
 Then there exists $n \leq \pent m2$ such that
$$0 \leq \alpha_0 \leq \alpha_1 \le \cdots \leq \alpha_n, \quad
 \alpha_n  \geq \alpha_{n+1} \geq \cdots \geq \alpha_{\pent m2} \geq 0.$$
\end{conjecture}
\pn
If this conjecture was true, it would imply the following property of
Alexander polynomials:
\begin{conjecture}\label{con-con}
Let $P(t) = a_0 -a_1 ( t+ t^{-1}) + a_2 ( t^2 + t^{-2} ) - \cdots
+(-1)^n a_n ( t^n + t^{-n} )$ be the Alexander polynomial of a
two-bridge knot. Then there exists an integer $k \le n$ such that 
$ (a_0, \ldots, a_k) $ is convex and $(a_{k}, \ldots, a_n)$ is concave.
\end{conjecture}
This property detects many non two-bridged polynomials
which are not detected by the other conditions.

\bibliographystyle{line}

\vfill
\pn
\hrule width 7cm height 1pt 
\pn
{\small
Pierre-Vincent {\sc Koseleff}\\
Universit{\'e} Pierre et Marie Curie (UPMC Paris 6) \& Inria-Rocquencourt\\
4, place Jussieu, F-75252 Paris Cedex 05 \\
e-mail: {\tt koseleff@math.jussieu.fr}
\pn
Daniel {\sc Pecker}\\
Universit{\'e} Pierre et Marie Curie (UPMC Paris 6),\\
4, place Jussieu, F-75252 Paris Cedex 05 \\
e-mail: {\tt pecker@math.jussieu.fr}}
\end{document}